\documentclass[12pt]{article}%
\usepackage{amsmath,amssymb}
\usepackage[applemac]{inputenc}
\usepackage{amsmath,amssymb,fullpage}
\usepackage{color}
\usepackage{amsmath}
\usepackage{amsfonts}
\usepackage{amssymb}
\usepackage{graphicx}%
\providecommand{\U}[1]{\protect\rule{.1in}{.1in}}

\newtheorem{definition}{Definition}[section]
\newtheorem{example}{Example}[section]
\newtheorem{theorem}[definition]{Theorem}
\newtheorem{problem}[definition]{Problem}
\newtheorem{remark}[definition]{ \it Remark}

\numberwithin{equation}{section}

\def\1B{\text{1\!\!I}}

\begin{document}

\date{27 April 2021}
\title{Optimal stopping, randomized stopping and singular control with general
information flow}
\author{Nacira Agram$^{1}$, Sven Haadem$^{2}$, Bernt \O ksendal$^{2}$ and Frank Proske$^{2}$}
\maketitle

\begin{abstract}
The purpose of this paper is two-fold:

\begin{itemize}
\item We extend the well-known relation between optimal stopping and
randomized stopping of a given stochastic process to a situation where the
available information flow is a filtration with no a priori assumed relation
to the filtration of the process. We call these problems optimal stopping and
randomized stopping with \emph{general} information.
\item Following an idea of Krylov \cite{K} we introduce a special
\emph{singular} stochastic control problem with general information and show
that this is also equivalent to the partial information optimal stopping and
randomized stopping problems. Then we show that the solution of this singular
control problem can be expressed in terms of partial information variational
inequalities.

\end{itemize}
\end{abstract}


\footnotetext[1]{Department of Mathematics, Linnaeus University, 35195 V\" axj\" o, Sweden. \newline
Email: \texttt{nacira.agram@lnu.se}}

\footnotetext[2]{Department of Mathematics, University of Oslo,
Blindern, 0316 Oslo, Norway. \newline Emails: \texttt{svenhaadem@gmail.com, oksendal@math.uio.no, proske@math.uio.no.}}



\section{Introduction}

There are several classic papers in the literature on the relation between
optimal stopping, randomized stopping and singular control of a given
stochastic process with filtration $\mathbb{F}:=\{\mathcal{F}_t\}_{t\geq0}$. See
e.g. Gy\"{o}ngy and
\v{S}i\v{s}ka \cite{GS} and the references therein. For other papers of
related interest, see Benth \& Reikvam \cite{BR}, Chaleyat-Maurel \textit{et al}
\cite{CMEKM}, Hu and \O ksendal \cite{HO}, Karatzas and Shreve \cite{KS}, \O ksendal \textit{et al}
\cite{O,OS}, and Shashiashvili \cite{Sha}.
Two fundamental references to the theory of optimal stopping are Dynkin \cite{Dyn} and Shiryaev \cite{Shi}.\\

The purpose of this paper is to extend this relation to a situation where the
admissible stopping times are required to be stopping times with respect to
another given \emph{information flow} $\mathbb{H}=\{\mathcal{H}_{t}\}_{t\geq
0}.$ We make no assumptions a priori about the relation between $\mathbb{H}$
and $\mathbb{F}$. If $\mathcal{H}_{t}\subseteq\mathcal{F}_{t}$ for all $t$, we
call this a \emph{partial information optimal stopping problem}.\newline If,
on the other hand $\mathcal{H}_{t}\supseteq\mathcal{F}_{t}$ for all $t$, we
call this an \emph{inside information optimal stopping problem}. Partial
information optimal stopping problems are studied in \O ksendal and Sulem \cite{OS2},
using a maximum principle for singular stochastic control of jump diffusions
and associated reflected backward differential equations. A special inside
information optimal stopping problem is studied (and solved) in Hu and \O ksendal \cite{HO}, based on Malliavin calculus and forward integration theory.

In the current paper, the admissible controls, such as the singular and the optimal stopping controls, are required to be $\mathbb{H}$-adapted. This is a
common situation in many applications, and one of our motivations for this
paper is to be able to study such more realistic optimal stopping problems.
In the current paper we extend the results of Gy\"{o}ngy and
\v{S}i\v{s}ka \cite{GS} and \O ksendal and Sulem \cite{OS2} to a more general
setting.
More precisely, we prove the equivalence of the following 3 problems: 
\begin{itemize}
\item Optimal stopping with general information flow.
\item Randomized stopping with general information flow.
\item Singular control with general information flow.
\end{itemize}
We then illustrate our result by finding explicit optimal stopping control under delayed information.
Finally, we obtain variational inequalities for singular control under partial information.


\section{Framework and problem formulations}

\vskip 0.5cm

Let $(\Omega,\mathcal{F},\mathbb{F}=\{\mathcal{F}_{t}\}_{t\geq0},P)$ be a
filtered probability space satisfying the usual conditions. Let $T\leq\infty$
be a fixed terminal time and let $\mathbb{H}:=\{\mathcal{H}_{t}\}_{t\geq0}$ be
another collection of complete $\sigma$-algebras $\mathcal{H}_{t}$,
not necessarily satisfying the usual conditions. 

\noindent We do not assume a priori that there is any relation between
$\mathbb{H}$ and $\mathbb{F}.$\newline For example, we could have

\begin{itemize}
\item $\mathcal{H}_{t}=\mathcal{F}_{(t-\delta)^{+}}$; $t\geq0$
(delayed/partial information case), or

\item $\mathcal{H}_{t}=\mathcal{F} _{t+\delta}$; $t\geq0$ with $\delta> 0$
(advanced information case).
\end{itemize}

\noindent Further, let $\mathcal{T}_{\mathbb{H}}=\mathcal{T}_{\mathbb{H}%
}^{(T)}$ denote the set of all $\mathbb{H}$-stopping times $\tau\leq T$, i.e.
the set of all functions
\[
\tau:\Omega\rightarrow\lbrack0,T],
\]
such that $\{\omega:\tau(\omega)\leq t\}\in\mathcal{H}_{t}$ for all
$t\in\lbrack0,T]$. In the following we let $\{k(t)\}_{t\geq0}$ be a given
$\mathbb{F}$-predictable process which is continuous at $t=0$. We assume that
$t\mapsto k(t)$ is 
left continuous with right sided limits  for all
$t\in\lbrack0,T]$ (c\` agl\` ad)
and satisfies
\begin{align}\label{2.1}
\underset{\tau\in\mathcal{T}_{\mathbb{H}}}{\sup}E[\lvert k(\tau)|]=:\kappa
<\infty.
\end{align}

\noindent We put $k(\tau(\omega))=0$ if $\tau(\omega)=\infty.$

\begin{remark}
If the filtration $\mathbb{H}$ satisfies the usual conditions, one can reduce
the problem to the complete information case when $\mathbb{F}=\mathbb{H}$ by
replacing the process $k(t)$ by its $\mathbb{H}$-optional conditional
expectation $\tilde{k}(t):=E[k(t) | \mathcal{H}_{t}]$. However, if $\mathbb{H}$ is a
strict subset of $\mathbb{F}$, we cannot go the other way.  More
precisely, given two arbitrary filtrations $\mathbb{F}$ and $\mathbb{H}$,
there is more information in the statement
\begin{equation}
\label{eq2.1}\underset{\tau\in\mathcal{T}_{\mathbb{H}}}{\sup} E[k(\tau)] =
\underset{G\in\mathcal{G}_{\mathbb{H}}}{\sup}E\left[  \int_{0} ^{T}%
k(t)dG(t)\right]  \text{ for any } \mathbb{F}\text{-adapted process } k(\cdot)
\end{equation}
than in the statement
\begin{equation}
\label{eq2.2}\underset{\tau\in\mathcal{T}_{\mathbb{H}}}{\sup} E[k(\tau)] =
\underset{G\in\mathcal{G}_{\mathbb{H}}}{\sup}E\left[  \int_{0}^{T}%
k(t)dG(t)\right]  \text{ for all } \mathbb{H}\text{-adapted process }
k(\cdot).
\end{equation}
Moreover, such a reduction may not be an advantage when it comes to solving
the problem. See Example 4.1.\newline We also point out that
several of our results do not need that the filtration $\mathbb{H}$ satisfies
the so-called "usual conditions", which would be needed for the reduction
argument above.
\end{remark}

 \emph{Note: All integrals in this paper are interpreted in the
Lebesgue-Stieltjes sense.}
\vskip 0.5cm
The purpose of the current paper is to study the relation between the following 3 problems in a general information flow context:\newline
We first consider the following \emph{general} information optimal stopping problem:

\begin{problem}{\bf(Optimal stopping)}
\label{problem1}\text{}\newline Find $\Phi\in\mathbb{R}$ and $\tau^{*}
\in\mathcal{T}_{\mathbb{H}}$ such that
\begin{align}
\Phi:= \underset{\tau\in\mathcal{T}_{\mathbb{H}}}{\sup} E[k(\tau)] =
E[k(\tau^{*})].
\end{align}

\end{problem}

Next we formulate the corresponding \emph{general} information randomized
stopping problem:

\begin{problem}{\bf(Randomized stopping)}\newline
\label{problem2} Let $\mathcal{G}_{\mathbb{H}}$ be the set of $\mathbb{H}%
$-adapted,  right-continuous and non-decreasing
processes $G(t); t \in[0,T]$ such that
\[
G(0)=0\text{ and } G(T) \leq 1\text{ a.s.}
\]
Find $\Lambda\in\mathbb{R}$ and $G^{\ast}\in\mathcal{G}_{\mathbb{H}}$ such
that
\[
\Lambda:=\underset{G\in\mathcal{G}_{\mathbb{H}}}{\sup}E\left[  \int_{0}%
^{T}k(t)dG(t)\right]  =E\left[  \int_{0}^{T}k(t)dG^{\ast}(t)\right]  .
\]

\end{problem}

\noindent Finally, we introduce our corresponding \emph{general} information
singular control problem:

\begin{problem}{\bf(Singular control)}\newline
\label{problem3} Let $\mathcal{A}_{\mathbb{H}}$ denote the set of all
$\mathbb{H}$-adapted
non-decreasing right-continuous processes
$\xi(t):[0,T]\rightarrow\lbrack0,\infty)$ such that $\xi(0)=0$ and
\[
\int_{[0,T]} \exp(-\xi(s))d\xi(s) \leq 1.
\]
 Find $\Psi\in\mathbb{R}$ and $\xi^{\ast}\in\mathcal{A}%
_{\mathbb{H}}$ such that
\[
\Psi:=\underset{\xi\in\mathcal{A}_{\mathbb{H}}}{\sup}E\left[  \int_{0}%
^{T}k(t)\exp\left(  -\xi(t)\right)  d\xi(t)\right]  = E\left[  \int_{0}%
^{T}k(t)\exp\left(  -\xi^{\ast}(t)\right)  d\xi^{\ast}(t)\right]  .
\]

\end{problem}

We will prove that all these 3 problems are equivalent, in the sense that
\[
\Phi=\Lambda=\Psi,
\]
and we will find explicit relations between the optimal $\tau^{\ast}$,
$G^{\ast}$and $\xi^{\ast}$.

\vskip 0.5cm

\section{Randomized stopping and optimal stopping with general information
flow}

In this section we prove that Problem \ref{problem1} and Problem
\ref{problem2} are equivalent. The following result may be regarded as an
extension of Theorem 2.1 in Gy\"{o}ngy and \v{S}i\v{s}ka \cite{GS} to general information:

\begin{theorem}
\label{theorem3.1}\text{{}}\newline%
\[
\Lambda:=\underset{G\in\mathcal{G}_{\mathbb{H}}}{\sup}E\left[  \int_{0}%
^{T}k(t)dG(t)\right]  =\underset{\tau\in\mathcal{T}_{\mathbb{H}}}{\sup
}E\left[  k(\tau)\right]  =:\Phi.
\]

\end{theorem}

\noindent{Proof.} \quad\text{ }\newline Choose $\tau\in\mathcal{T}%
_{\mathbb{H}}$ and define, for $n = 1,2, \ldots$,
\begin{align}%
\begin{cases}
G^{(n)}(t) & =\mathbf{1}_{\{t\geq\tau>0\}}+(1-e^{-nt})\mathbf{1}_{\{\tau
=0\}},\text{ for } t < T,\\
G^{(n)}(T) & = 1.
\end{cases}
\end{align}

Then $G^{(n)}(\cdot)\in\mathcal{G}_{\mathbb{H}}$ and we see that
\[
E\left[  k(\tau)\right]  =\underset{n\rightarrow\infty}{\lim}E\left[  \int
_{0}^{T}k(t)dG^{(n)}(t)\right]  \leq\underset{G\in\mathcal{G}_{\mathbb{H}}%
}{\sup}E\left[  \int_{0}^{T}k(t)dG(t)\right]  .
\]
Since $\tau\in\mathcal{T}_{\mathcal{H}}$ was arbitrary, this proves that
\[
\underset{\tau\in\mathcal{T}_{\mathbb{H}}}{\sup}E\left[  k(\tau)\right]
\leq\underset{G\in\mathcal{G}_{\mathbb{H}}}{\sup}E\left[  \int_{0}%
^{T}k(t)dG(t)\right]  .
\]
To get the opposite inequality, we define for each $G\in\mathcal{G}%
_{\mathbb{H}}$ and $r\in\lbrack0,G(T))=[0,1)$, the time change $\alpha(r)$ by
\[
\alpha(r)=\inf\{s\geq0;G(s) \geq r\}.
\]
Then $\{\omega;\alpha(r) \leq t\}=\{\omega;G(t) \geq r\}\in\mathcal{H}_{t}$,
so $\alpha(r)\in\mathcal{T}_{\mathbb{H}}$ for all $r$. Moreover,
$G(\alpha(t))=t$ for a.a. $t$ and hence
\[
E\left[  \int_{0}^{T}k(t)dG(t)\right]  =E\left[  \int_{0}^{G(T)}%
k(\alpha(r))dr\right]  \leq\int_{0}^{1}\underset{\tau\in\mathcal{T}%
_{\mathbb{H}}}{\sup}E[k(\tau)]dr=\underset{\tau\in\mathcal{T}_{\mathbb{H}}%
}{\sup}E\left[  k(\tau)\right]  .
\]
\hfill$\square$ \bigskip\vskip0.5cm

\section{Singular control and optimal stopping with general
information}

In this section we prove that Problem \ref{problem1} and Problem
\ref{problem3} are equivalent:

\begin{theorem}
\label{theorem4.2} 
Define $\mathcal{A}_{\mathbb{H}}^{c}=\{\xi \in \mathcal{A}_{\mathbb{H}};\xi $
is continuous$\}$, $\overline{\mathbb{H}}=\{\overline{\mathcal{H}}%
_{t}\}_{0\leq t\leq T}=\{\mathcal{H}_{t}\cap \mathcal{F}_{t^{-}}\}_{0\leq
t\leq T}$ ($\mathcal{F}_{0^{-}}:=\mathcal{F}_{0}$), $\mathcal{G}_{\overline{%
\mathbb{H}}}^{\ast }=\{G\in \mathcal{G}_{\overline{\mathbb{H}}};G$ is $%
\overline{\mathbb{H}}-$predictable$\}$ and $\mathcal{T}_{\overline{\mathbb{H}%
}}^{\ast }=\{\tau \in \mathcal{T}_{\overline{\mathbb{H}}};\tau $ is $%
\overline{\mathbb{H}}-$predictable$\}$. Further, we assume that the
information flow $\overline{\mathbb{H}}$ is right continuous and that $%
E[\sup_{0\leq t\leq T}\left\vert k(t)\right\vert ]<\infty $. Then%
\begin{eqnarray*}
\sup_{\xi \in \mathcal{A}_{\mathbb{H}}^{c}}E\Big[\int_{0}^{T}k(t)\exp (-\xi
(t))d\xi (t)\Big] &=&\Psi :=\sup_{\xi \in \mathcal{A}_{\mathbb{H}%
}}E\Big[\int_{0}^{T}k(t)\exp (-\xi (t))d\xi (t)\Big]\\
&=&\sup_{G\in \mathcal{G}_{\mathbb{H}}}E\Big[\int_{0}^{T}k(t)dG(t)\Big] =\sup_{\tau
\in \mathcal{T}_{\mathbb{H}}}E[k(\tau )] \\
&=:&\Phi= \sup_{G\in \mathcal{G}_{\overline{\mathbb{H}}}^{\ast
}}E\Big[\int_{0}^{T}k(t)dG(t)\Big] =\sup_{\tau \in \mathcal{T}_{\overline{\mathbb{H%
}}}^{\ast }}E[k(\tau )].
\end{eqnarray*}
\end{theorem}

 \noindent{Proof.} \quad\text{ }\newline Let $\xi\in
\mathcal{A}_{\mathbb{H}}$.  Then $w(t):=\int_{[0,t]} e^{-\xi
(s)}d\xi(s) \in\mathcal{G}_{\mathbb{H}}$ and hence, by
Theorem \ref{theorem3.1},
\begin{align*}
E\left[  \int_{0}^{T}k(t)\exp\left(  -\xi(t)\right)  d\xi(t)\right]   &
=E\left[  \int_{0}^{T}k(t)dw(t)\right]  \leq\underset{G\in\mathcal{G}%
_{\mathbb{H}}}{\sup}E\left[  \int_{0}^{T}k(t)dG(t)\right] \\
&  =\underset{\tau\in\mathcal{T}_{\mathbb{H}}}{\sup}E\left[  k(\tau)\right]  .
\end{align*}
Therefore,
\begin{equation}
\underset{\xi\in\mathcal{A}_{\mathbb{H}}}{\sup}E\left[  \int_{0}^{T}%
k(t)\exp(-\xi(t))d\xi(t)\right]  \leq\underset{\tau\in\mathcal{T}_{\mathbb{H}%
}}{\sup}E\left[  k(\tau)\right]  . \label{eq:supV_pos_k}%
\end{equation}
To get the opposite inequality, choose $\delta > 0$, and $\tau\in\mathcal{T}_{\mathbb{H}}$. By left-continuity of $k$ we may assume that $\tau < T$ a.s.
Define, for $n=1,2,\ldots$
\begin{align}
u^{(n)}(t)=%
\begin{cases}
0 & \text{ for }t<\tau,\\
n & \text{ for }t\geq\tau,
\end{cases}
\end{align}
and
\begin{align}
G^{(n)}(t)=%
\begin{cases}
0 & \text{ for }t<\tau,\\
1-e^{-n(t-\tau)} & \text{ for }t\geq \tau.
\end{cases}
\end{align}
Then $\xi^{(n)}(t):=\int_{0}^{t}u^{(n)}(s)ds\in\mathcal{A}_{\mathbb{H}}^{c} $,
$G^{(n)}(t)\in\mathcal{G}_{\mathbb{H}}$ and for any $\delta>0$, we have
 \begin{align*}
&\int_{0}^{T}k(t)u^{(n)}(t)\exp\left(  -\int_{0}^{t}u^{(n)}(s)ds\right)dt\\
&=\int_{0}^{T}k(t^{+})u^{(n)}(t)\exp\left(  -\int_{0}^{t}u^{(n)}(s)ds\right)dt\\
&=\int_{0}^{T}k(t^{+})dG^{(n)}(t)=\int_{\tau}^{T}k(t^{+})dG^{(n)}(t)
=I_{n}+J_{n}+K_{n},
\end{align*}
where
\begin{align}
I_{n}&=\int_{\tau}^{(\tau+ \delta) \wedge T}k(\tau^{+})dG^{(n)}(t),\nonumber\\
J_{n}&=\int_{\tau}^{(\tau+ \delta) \wedge T}\left(  k(t^{+})-k(\tau^{+})\right)  dG^{(n)}(t),\nonumber\\
K_{n}&=\int_{(\tau+ \delta) \wedge T}^{T}k(t^{+})dG^{(n)}(t),
\end{align}
so we see that when $\tau+ \delta<T$
\[
I_{n}=k(\tau^{+})(1-e^{-n\delta})\rightarrow
k(\tau^{+})\text{ as }n\rightarrow\infty,
\]%
and when $\tau+ \delta>T$
\[
I_{n}=k(\tau^{+})(1-e^{-n(T-\tau)}) \rightarrow k(\tau^{+})\text{ as }n\rightarrow\infty.
\]
By right-continuity,
\[
\lvert J_{n}\rvert\leq\sup_{t\in\lbrack\tau,\tau+\delta]}\lvert k(t^{+})-k(\tau^{+}
)\rvert\rightarrow0\text{ when }\delta\rightarrow0.
\]
Moreover, 
\begin{align*}
|K_{n}|   &  \leq\sup_{t\in [0,T]}
|k(t^{+})|  [e^{-n((\tau + \delta)\wedge T-\tau)}-e^{-n(T-\tau)} ]\rightarrow 0\text{ when }n\rightarrow\infty.
\end{align*}
Combining the above in connection with $E[\sup_{0\leq t\leq T}\left\vert
k(t)\right\vert ]<\infty $, it follows from dominated convergence theorem, that%
\begin{equation*}
\lim_{\delta \longrightarrow 0}\lim_{n\longrightarrow \infty
}E\Big[\int_{0}^{T}k(t)u^{(n)}(t)\exp (-\int_{0}^{t}u^{(n)}(s)ds)dt\Big]=E\Big[k(\tau
^{+})\Big].
\end{equation*}%
Therefore,%
\begin{equation*}
\sup_{\xi \in \mathcal{A}_{\mathbb{H}}^{c}}E\Big[\int_{0}^{T}k(t)\exp (-\xi
(t))d\xi (t)\Big]\geq E[k(\tau ^{+})].
\end{equation*}%
Since $\tau \in \mathcal{T}_{\mathbb{H}}$ was arbitrary this proves that%
\begin{equation}
\sup_{\xi \in \mathcal{A}_{\mathbb{H}}^{c}}E[\int_{0}^{T}k(t)\exp (-\xi
(t))d\xi (t)]\geq \sup_{\tau \in \mathcal{T}_{\mathbb{H}}}E[k(\tau ^{+})].
\label{InEq-1}
\end{equation}%
Let $\tau \in \mathcal{T}_{\overline{\mathbb{H}}}^{\ast }$. Then, using the
fact that%
\begin{equation*}
\lim_{\substack{ t\longrightarrow t_{0} \\ t<t_{0}}}k(t^{+})=k(t_{0}),t_{0}%
\in (0,T],
\end{equation*}%
we can find by assumption an announcing sequence of stopping times $\left(
\tau _{n}\right) _{n\geq 1}\subset \mathcal{T}_{\overline{\mathbb{H}}}$ such
that $\tau _{n}\leq \tau $ increases to $\tau $ and $\tau _{n}<\tau ,$
whenever $\tau >0$. So $k(\tau _{n}^{+})\underset{n\longrightarrow \infty }{%
\longrightarrow }k(\tau )$ a.e. Then dominated convergence theorem, yields%
\begin{equation*}
E[k(\tau _{n}^{+})]\underset{n\longrightarrow \infty }{\longrightarrow }%
E[k(\tau )].
\end{equation*}%
So%
\begin{equation*}
\sup_{\tau \in \mathcal{T}_{\mathbb{H}}}E[k(\tau ^{+})]\geq E[k(\tau ^{\ast
})],\tau ^{\ast }\in \mathcal{T}_{\overline{\mathbb{H}}}^{\ast },
\end{equation*}%
which implies%
\begin{equation}
\sup_{\tau \in \mathcal{T}_{\mathbb{H}}}E[k(\tau ^{+})]\geq \sup_{\tau \in 
\mathcal{T}_{\overline{\mathbb{H}}}^{\ast }}E[k(\tau )].  \label{InEq0}
\end{equation}%
If $G\in \mathcal{G}_{\overline{\mathbb{H}}}^{\ast }$, we see from the proof
of Theorem \ref{theorem3.1} that%
\begin{equation*}
\alpha (r)\in \mathcal{T}_{\overline{\mathbb{H}}}^{\ast },r\in \lbrack 0,1).
\end{equation*}%
So 
\begin{equation}
\sup_{G\in \mathcal{G}_{\overline{\mathbb{H}}}^{\ast
}}E\Big[\int_{0}^{T}k(t)dG(t)\Big]\leq \sup_{\tau \in \mathcal{T}_{\overline{\mathbb{%
H}}}^{\ast }}E\Big[k(\tau )\Big]  \label{InEq1}
\end{equation}%
by the proof of Theorem \ref{theorem3.1}. Further, we observe that%
\begin{equation}
\sup_{\xi \in \mathcal{A}_{\mathbb{H}}^{c}}E\Big[\int_{0}^{T}k(t)\exp (-\xi
(t))d\xi (t)\Big]\leq \sup_{G\in \mathcal{G}_{\mathbb{H}}}E\Big[%
\int_{0}^{T}k(t)dG(t)\Big].  \label{InEq2}
\end{equation}%
On the other hand, we know in connection with our assumptions that%
\begin{eqnarray*}
E\Big[\int_{0}^{T}k(t)dG(t)\Big] &=&E\Big[\int_{0}^{T}\left. ^{p}(k)(t)\right. dG(t)\Big] 
=E\Big[\int_{0}^{T}k(t)d(G)^{p}(t)\Big],
\end{eqnarray*}%
where $^{p}(\cdot )$ and $(\cdot )^{p}$ denote the predictable and the dual
predictable projection with respect to the filtration $\mathbb{F}$,
respectively. Under our assumptions on $G$ it is known that $%
(G)^{p}(t),0\leq t\leq T$ is right continuous and non-decreasing (see
Dellacherie, Meyer \cite{DM}). Further, we see from the definition of $(\cdot )^{p}$
that $(G)^{p}(t),0\leq t\leq T$ is $\overline{\mathbb{H}}-$adapted with $%
(G)^{p}(0)=0$ and $(G)^{p}(T)\leq 1$ a.e. On the other hand, since $%
\overline{\mathbb{H}}$ is right continuous, we find that elementary $\mathbb{%
F}-$predictable processes, which are $\overline{\mathbb{H}}-$adapted, are $%
\overline{\mathbb{H}}-$predictable. Using the latter fact combined with the
monotone class theorem, it follows that $(G)^{p}(t),0\leq t\leq T$ belongs
to $\mathcal{G}_{\overline{\mathbb{H}}}^{\ast }$. So%
\begin{equation}
\sup_{G\in \mathcal{G}_{\mathbb{H}}}E\Big[\int_{0}^{T}k(t)dG(t)\Big]=\sup_{G\in \mathcal{G}_{\mathbb{H}}}E\Big[\int_{0}^{T}k(t)d(G)^p(t)\Big]=\sup_{G\in 
\mathcal{G}_{\overline{\mathbb{H}}}^{\ast }}E\Big[\int_{0}^{T}k(t)dG(t)\Big].
\label{Eq1}
\end{equation}%
Therefore, (\ref{InEq1}) and (\ref{InEq2}) entail that%
\begin{eqnarray}\label{4.10a}
\sup_{\xi \in \mathcal{A}_{\mathbb{H}}^{c}}E\Big[\int_{0}^{T}k(t)\exp (-\xi
(t))d\xi (t)\Big] &\leq &\sup_{G\in \mathcal{G}_{\overline{\mathbb{H}}}^{\ast
}}E\Big[\int_{0}^{T}k(t)dG(t)\Big] \nonumber\\
&\leq &\sup_{\tau \in \mathcal{T}_{\overline{\mathbb{H}}}^{\ast }}E[k(\tau
)].
\end{eqnarray}

So we conclude from (\ref{InEq-1}), (\ref{eq:supV_pos_k}), Theorem \ref{theorem3.1}, (\ref{Eq1}) and \eqref{4.10a} and (\ref{InEq0}) (in that order) that
\begin{align*}
&\sup_{\tau \in \mathcal{T}_{\mathbb{H}}}E[k(\tau^+)]\\
\text{ (by (\ref{InEq-1}) )} \hskip 0.5cm &\leq \sup_{\xi \in \mathcal{A}_{\mathbb{H}}^{c}}E\Big[\int_{0}^{T}k(t)\exp (-\xi
(t))d\xi (t)\Big] \\
( \text{ since } \mathcal{A}_{\mathbb{H}}^c \subset \mathcal{A}_{\mathbb{H}}  ) \hskip 0.5cm  &\leq \sup_{\xi \in \mathcal{A}_{\mathbb{H}}}E\Big[\int_{0}^{T}k(t)\exp (-\xi
(t))d\xi (t)\Big] \\
(\text{ by } (\ref{eq:supV_pos_k}) ) \hskip 0.5cm &\leq \sup_{\tau \in \mathcal{T}_{\mathbb{H}}}E[k(\tau )] \\
(\text{ by  Theorem } \ref{theorem3.1} ) \hskip 0.5cm &=\sup_{G\in \mathcal{G}_{\mathbb{H}}}E\Big[\int_{0}^{T}k(t)dG(t)\Big]\\
( \text{ by } (\ref{Eq1}) ) \hskip 0.5cm &=\sup_{G\in \mathcal{G}_{\overline{\mathbb{H}}}^{\ast }}E\Big[\int_{0}^{T}k(t)dG(t)\Big]\\
(\text{by } \eqref{4.10a}) \hskip 0.5cm &\leq \sup_{\tau \in \mathcal{T}_{\overline{\mathbb{H}}}^{\ast }}E[k(\tau)]\\
(\text{ by } (\ref{InEq0}) \hskip 0.5cm &\leq \sup_{\tau \in \mathcal{T}_{\mathbb{H}}}E|k(\tau^+)].
\end{align*}
Since the first term in this chain of inequalities/equations is the same as the last term, we conclude that all the terms are the same and the proof follows. 

\hfill$\square$ \bigskip

\noindent It is of interest to find the connection between an optimal stopping
time $\tau^{\ast}\in\mathcal{T}_{\mathbb{H}}$ for Problem \ref{problem1} and
the corresponding optimal singular controls $G^{\ast}$, $\xi^{\ast}$ for
Problem \ref{problem2} and Problem \ref{problem3}, respectively. The
connection is given by the following result:

\begin{theorem}
\begin{description}
\item[a)] Suppose $\tau^{\ast}\in\mathcal{T}_{\mathbb{H}}$ is an optimal
stopping time for Problem \ref{problem1}. Define
\begin{equation}
G^{\ast}(t):=\mathbf{1}_{\{t\geq\tau^{\ast}>0\}}+\mathbf{1}_{\{\tau^{\ast
}=0\}}. \label{eq:G_pos_k}%
\end{equation}
Then $G^{\ast}\in\mathcal{G}_{\mathbb{H}}$ is an optimal singular control for
Problem \ref{problem2}. \vskip0.3cm

\item[b)] Conversely, suppose $G^{\ast}\in\mathcal{G}_{\mathbb{H}}$ is an
optimal singular control for Problem \ref{problem2}. Define 
\begin{equation}
\alpha^{\ast}(r):=\inf\{s\geq0;G^{\ast}(s) \geq r\};\text{ for }r\in
\lbrack0,1). \label{eq:alpha_pos_k}%
\end{equation}
\color{black} Then $\alpha^{\ast}(r)\in\mathcal{T}_{\mathbb{H}}$ and
$\alpha^{\ast}(r) $ is an optimal stopping time for Problem \ref{problem1},
for all $r\in\lbrack0,1)$.

\item[c)] Suppose $\xi^{\ast}\in\mathcal{A}_{\mathbb{H}}$ is an optimal
control for Problem \ref{problem3}.  Then the process
\[
G^{\ast}(t):=\int_{[0,t]}\exp(-\xi^{\ast}(s))d\xi^{\ast}(s)
\]
is an optimal control for Problem \ref{problem2}. \color{black}

\item[d)] Conversely, suppose $G^{\ast}(t)$ is an optimal control for Problem
\ref{problem2}. Define $\xi^{\ast}(t)$ to be a solution of the differential
equation 
\[
d\xi^{\ast}(t)=\exp(\xi^{\ast}(t))dG^{\ast}(t);\text{ }\xi^{\ast}(0^{-})=0.
\]
\color{black} Then $\xi^{\ast}(t)$ is an optimal control for Problem
\ref{problem3}.
\end{description}
\end{theorem}

\vskip 0.3cm

\noindent{Proof.}

\begin{enumerate}
\item[a)] Suppose $\tau^{\ast}\in\mathcal{T}_{\mathbb{H}}$ is optimal for
Problem \ref{problem1} and let $G^{\ast}$ be as in \eqref{eq:G_pos_k}. Then by
Theorem \ref{theorem3.1}
\begin{align*}
\sup_{\tau\in\mathcal{T}_{\mathbb{H}}}E\left[  k(\tau)\right]   &  =E\left[
k(\tau^{\ast})\right]  =E\left[  \int_{0}^{T}k(t)dG^{\ast}(t)\right] \\
&  \leq\sup_{G\in\mathcal{G}_{\mathbb{H}}}E\left[  \int_{0}^{T}%
k(t)dG(t)\right]  =\sup_{\tau\in\mathcal{T}_{\mathbb{H}}}E\left[
k(\tau)\right]  .
\end{align*}
Hence we have equality in the above, and therefore
\[
E\left[  \int_{0}^{T}k(t)dG^{\ast}(t)\right]  =\sup_{G\in\mathcal{G}%
_{\mathbb{H}}}E\left[  \int_{0}^{T}k(t)dG(t)\right]  ,
\]
which proves that $G^{\ast}$ is optimal for Problem \ref{problem2}.\newline

\item[b)] Conversely, suppose $G^{\ast}\in\mathcal{G}_{\mathbb{H}}$ is optimal
for Problem \ref{problem2}. Let $\alpha^{\ast}(r)$ be as in
\eqref{eq:alpha_pos_k}. Then $\alpha^{\ast}(r)\in\mathcal{T}_{\mathbb{H}}$ for
all $r$ and, by Theorem \ref{theorem3.1},
\begin{align*}
\sup_{G\in\mathcal{G}_{\mathbb{H}}}E\left[  \int_{0}^{T}k(t)dG(t)\right]   &
=E\left[  \int_{0}^{T}k(t)dG^{\ast}(t)\right] \\
&  =E\left[  \int_{0}^{G^{\ast}(T)}k(\alpha^{\ast}(r))dr\right]  =\int_{0}%
^{1}E\left[  k(\alpha^{\ast}(r))\right]  dr\\
&  \leq\int_{0}^{1}\sup_{\tau\in\mathcal{T}_{\mathbb{H}}}E\left[
k(\tau)\right]  dr=\sup_{\tau\in\mathcal{T}_{\mathbb{H}}}E\left[
k(\tau)\right] \\
&  \leq\sup_{G\in\mathcal{G}_{\mathbb{H}}}E\left[  \int_{0}^{T}%
k(t)dG(t)\right]  .
\end{align*}
We conclude that we have equality everywhere in the above, and therefore
\begin{equation}
\int_{0}^{1}E\left[  k(\alpha^{\ast}(r))\right]  dr=\sup_{\tau\in
\mathcal{T}_{\mathbb{H}}}E\left[  k(\tau)\right]  . \label{eq4.6a}%
\end{equation}
Since $\alpha^{\ast}(r)$ is a stopping time for all $r\in\lbrack0,1)$ we have
\begin{equation}
E\left[  k(\alpha^{\ast}(r))\right]  \leq\sup_{\tau\in\mathcal{T}_{\mathbb{H}%
}}E\left[  k(\tau)\right]  \quad\forall r.
\end{equation}
Therefore \eqref{eq4.6a} is only possible if
\begin{equation}
E\left[  k(\alpha^{\ast}(r))\right]  =\sup_{\tau\in\mathcal{T}_{\mathbb{H}}%
}E\left[  k(\tau)\right]  ,\text{ for a.a. }r\in\lbrack0,1). \label{eq4.6}%
\end{equation}

Choose arbitrary $\bar{r}\in (0,1]$. Then since $\alpha^{\ast}(r)$ is
left-continuous 
we can find $r_{n}\in(0,1)$ such that $\alpha^{\ast}(r_{n})$
is optimal for all $n$ and $\alpha^{\ast}(r_{n})\rightarrow\alpha^{\ast}%
(\bar{r})$ as $n\rightarrow\infty$. This gives
\[
E\left[  k(\alpha^{\ast}(\bar{r}))\right]  =\lim_{n\rightarrow\infty}E\left[
k(\alpha^{\ast}(r_{n}))\right]  =\sup_{\tau\in\mathcal{T}_{\mathbb{H}}%
}E\left[  k(\tau)\right]  .
\]
\color{black} Hence $\alpha^{\ast}(r)$ is an optimal stopping time for all
$r\in (0,1]$.

\item[c),d)] If $G^{\ast}(t)$ and $\xi^{\ast}(t)$ are chosen as given in c)
and d) respectively, then we see in either case that
\[
E\Big[\int_{0}^{T}k(s)dG^{\ast}(s)\Big]=E\Big[\int_{0}^{T}k(s)\exp(-\xi^{\ast}%
(s))d\xi^{\ast}(s)\Big].
\]
The two statements c) and d) follow from this.
\end{enumerate}

\hfill$\square$ \bigskip

\begin{remark}
In the case when
\[
k(t)\geq0,\text{ for all }t\in\lbrack0,T],
\]
the optimal $G^{\ast}\in\mathcal{G}_{\mathbb{H}}$ satisfies $$G^{\ast}(T)=1,$$ and the optimal $\xi\in\mathcal{A}_{\mathbb{H}}$ satisfies $$\int_{\lbrack0,T]}\exp(-\xi^{\ast}(s))d\xi^{\ast}(t)=1.$$
\end{remark}

\begin{example}
{\bf{(The optimal time to sell when there is delayed information)}}\newline 
To illustrate our results, we follow the Example 3.1 in \O ksendal \cite{O}. \newline Let $\mathbb{F}$ be the filtration of
one-dimensional Brownian motion $B(\cdot)$ and let $\mathbb{H}=\{\mathcal{H}%
_{t}\}$ be the delayed information flow given by $\mathcal{H}_{t}%
=\mathcal{F}_{(t-\delta)^{+}}$ for some constant $\delta>0$. 
Define
\[
k(t)=e^{-\rho t}(X(t)-a);\quad t\in\lbrack0,T],
\]
where $\rho>0,a>0$ are given constants, and the process $X$ is a geometric
Brownian motion of the form
\[
dX(t)=X(t)[\mu(t)dt+\sigma(t)dB(t)];\quad X(0)>0,
\]
where $\mu(t)$ and $\sigma(t)>0$ are bounded $\mathbb{F}$-adapted processes.
Then it follows from Theorem 3.1 in \cite{O} that the optimal stopping time
$\tau^{\ast}$ for Problem 2.2 has the form
\[
\tau^{\ast}=\alpha+\delta,
\]
where $\alpha$ is the optimal stopping time (which in some cases can be found
explicitly) for a related optimal stopping problem with non-delayed
information. \color{black} By Theorem 4.3 a) the optimal $G^{\ast}$ for the
corresponding randomized stopping problem (Problem 2.3) is
\[%
\begin{cases}
G^{\ast}(t)=0\quad t<\tau^{\ast};\\
G^{\ast}(t)=1\quad t\geq\tau^{\ast}.
\end{cases}
\]
And by Theorem 4.3 d) the optimal $\xi^{\ast}$ for the corresponding singular
control problem (Problem 2.4) is given by 
\[
\exp(-\xi^{\ast}(t))d\xi^{\ast}(t)=\delta_{\tau^{\ast}}(t);\quad t\in
\lbrack0,T];
\]
where $\delta_{\tau^{\ast}}(t)$ is the Dirac point mass at $t=\tau^{\ast}$. 
\end{example}


\section{Singular control with partial information flow}
In this section we assume
that $\mathbb{H}$ satisfies the usual conditions and that we are in a partial information setting, i.e. that
\begin{equation}
\mathcal{H}_{t}\subseteq\mathcal{F}_{t}\text{ for all }t.
\end{equation}
We also suppose that the terminal time horizon $T < \infty$.
\\
We now turn to the partial information singular control problem
(Problem \ref{problem3}):

\begin{problem}
\label{problem5.1} Find $\Psi\in\mathbb{R}$ and $\hat{\xi}\in\mathcal{A}%
_{\mathbb{H}}$ such that
\begin{equation}
\Psi=\sup_{\xi\in\mathcal{A}_{\mathbb{H}}}J(\xi)=J(\hat{\xi}),
\label{eq:Psi2}%
\end{equation}
where
\begin{equation}
J(\xi)=E\left[  \int_{0}^{T}k(t)\exp\left(  -\xi(t)\right)  d\xi(t)\right]  .
\label{eq:J}%
\end{equation}
\end{problem}

\noindent Problem \ref{problem5.1} can be considered as a generalisation of
the singular control problem discussed in Section 2 of \O ksendal and Sulem
\cite{OS2}, where a singular control version of the maximum principle is used.
 However, since the singular control $\xi$ appears both in the
integrand and as an integrator, the problem \eqref{eq:Psi2} - \eqref{eq:J} is
not covered by the results in \cite{OS2}.
Here we give a direct approach based on a variational argument. \newline%
\newline For $\xi
\in\mathcal{A}_{\mathbb{H}}$, we define $\mathcal{V(\xi)}$ to be the set of
c\`adl\`ag processes $\zeta(t):[0,T]\rightarrow\lbrack0,\infty]$ of finite
variation such that there exists $\delta=\delta(\xi)>0$ such that
\[
\xi+y\zeta\in\mathcal{A}_{\mathbb{H}}\text{ for all }y\in\lbrack0,\delta].
\]
For $\xi\in\mathcal{A}_{\mathbb{H}}$ and $\zeta\in\mathcal{V}(\xi)$, we define
$D(\xi,\zeta)\in\mathbb{R}$ by%
\small
\begin{align}
&  D(\xi,\zeta):=\lim_{y\rightarrow0^{+}}\sup\frac{1}{y}\left(  J(\xi
+y\zeta)-J(\xi)\right) \nonumber\\
&  =\lim_{y\rightarrow0^{+}}\sup\frac{1}{y}\left(  E\left[  \int_{0}%
^{T}k(s)\left\{  \exp\Big(-\big(\xi(s)+y\zeta(s)\big)\Big)\Big(d\xi
(s)+yd\zeta(s)\Big)-\exp\big(-\xi(s)\big)d\xi(s)\right\}  \right]  \right)
\nonumber\\
&  =\lim_{y\rightarrow0^{+}}\sup\frac{1}{y}\Bigg(E\Bigg[\int_{0}%
^{T}k(s)\Big\{\exp\big(-\xi(s)\big)\big(\exp\left(  -y\zeta(s)\right)
-1\big)d\xi(s)\nonumber\\
&  +y\exp\big(-\xi(s)\big)\exp\big(-y\zeta(s)\big)d\zeta
(s)\Big\}\Bigg]\Bigg)\nonumber\\
&  =E\left[  \int_{0}^{T}k(s)\exp\big(-\xi(s)\big)\left\{  -\zeta
(s)d\xi(s)+d\zeta(s)\right\}  \right]  . \label{eq:xi,zeta}%
\end{align}
Now suppose $\xi=\hat{\xi}$ maximizes $J(\xi)$. Then by \eqref{eq:xi,zeta}
\begin{equation}
E\left[  \int_{0}^{T}k(s)\exp\left(  -\hat{\xi}(s)\right)  \left\{
-\zeta(s)d\hat{\xi}(s)+d\zeta(s)\right\}  \right]  =D(\hat{\xi},\zeta)\leq0,
\label{eq:xi,zeta2}%
\end{equation}
for all $\zeta\in\mathcal{V}(\hat{\xi})$. In particular, if we for $\delta>0$
choose
\[
\zeta_{0}(s)=%
\begin{cases}
0 & s<t,\\
\frac{(s-t)a}{\delta} & t\leq s\leq t+\delta,\\
a & s\geq t+\delta,
\end{cases}
\]
for some $t\in\lbrack0,T]$ and some bounded $\mathcal{H}_{t}$-measurable
random variable $a \geq 0$, then $\zeta_{0}\in
\mathcal{V}(\hat{\xi}) $ and \eqref{eq:xi,zeta2} gives
\[%
\begin{array}
[c]{c}%
E\left[  \int_{t}^{t+\delta}k(s)\exp\left(  -\hat{\xi}(s)\right)
\frac{(s-t)a}{\delta}d\hat{\xi}(s)+\int_{t+\delta}^{T}k(s)\exp\left(
-\hat{\xi}(s)\right)  ad\hat{\xi}(s)\right. \\
\left.  -\int_{t}^{t+\delta}k(s)\exp\left(  -\hat{\xi}(s)\right)  \frac
{a}{\delta}ds\right]  \geq0.
\end{array}
\]
Since this holds for all such $a$ and all $\delta>0$, we
conclude that
\[
E\left[  \int_{t}^{T}k(s)\exp\left(  -\hat{\xi}(s)\right)  d\hat{\xi
}(s)-k(t)\exp\left(  -\hat{\xi}(t)\right)  \bigg \lvert\mathcal{H}_{t}\right]
\geq0;\text{ }t\in\lbrack0,T].
\]
Next, let us choose

\begin{enumerate}
\item $d\zeta_{1}(s) = d\hat{\xi}(s)$ and

\item $d\zeta_{2}(s) = -d\hat{\xi}(s)$.
\end{enumerate}

\noindent Then $\zeta_{i}\in\mathcal{V}(\hat{\xi})$ for $i=1,2$ and
\eqref{eq:xi,zeta2} gives
\begin{equation}
E\left[  \int_{0}^{T}k(s)\exp\left(  -\hat{\xi}(s)\right)  \left\{  -\hat{\xi
}(s)d\hat{\xi}(s)+d\hat{\xi}(s)\right\}  \right]  =0. \label{eq:exp2}%
\end{equation}
Note that by the Fubini theorem we have
\begin{align}
&  \int_{0}^{T}\left(  \int_{t}^{T}k(s)\exp\left(  -\hat{\xi}(s)\right)
d\hat{\xi}(s)\right)  d\hat{\xi}(t)\label{eq:Fub}\\
&  =\int_{0}^{T}\left(  \int_{0}^{s}d\hat{\xi}(t)\right)  k(s)\exp\left(
-\hat{\xi}(s)\right)  d\hat{\xi}(s)\nonumber\\
&  =\int_{0}^{T}k(s)\exp\left(  -\hat{\xi}(s)\right)  \hat{\xi}(s)d\hat{\xi
}(s).\nonumber
\end{align}
Substituting \eqref{eq:Fub} into \eqref{eq:exp2} we get
\[
E\left[  \int_{0}^{T}\left\{  \int_{t}^{T}k(s)\exp\left(  -\hat{\xi
}(s)\right)  d\hat{\xi}(s)-k(t)\exp\left(  -\hat{\xi}(t)\right)  \right\}
d\hat{\xi}(t)\right]  =0.
\]
This proves part a) of the following theorem: \vskip0.5cm

\begin{theorem}
[Variational inequalities]\label{thm:Var_ineq}\text{ }\newline

\begin{enumerate}
\item[a)] Suppose $\hat{\xi}\in\mathcal{A}_{\mathbb{H}}$ is optimal for
\eqref{eq:Psi2} - \eqref{eq:J}. Then
\begin{equation}
E\left[  \int_{t}^{T}k(s)\exp\left(  -\hat{\xi}(s)\right)  d\hat{\xi
}(s)-k(t)\exp\left(  -\hat{\xi}(t)\right)  \bigg\lvert\mathcal{H}_{t}\right]
\geq0;\text{ }t\in\lbrack0,T], \label{eq:thm5.1_a_1}%
\end{equation}
and
\begin{equation}
E\left[  \int_{t}^{T}k(s)\exp\left(  -\hat{\xi}(s)\right)  d\hat{\xi
}(s)-k(t)\exp\left(  -\hat{\xi}(t)\right)  \bigg\lvert\mathcal{H}_{t}\right]
d\hat{\xi}(t)=0;\text{ }t\in\lbrack0,T]. \label{eq:thm5.1_a_2}%
\end{equation}

\item[b)] Conversely, suppose \eqref{eq:thm5.1_a_1} - \eqref{eq:thm5.1_a_2}
hold for some $\hat{\xi}\in\mathcal{A}_{\mathbb{H}}$. Then
\begin{equation}
D(\hat{\xi},\zeta)\leq0\text{ for all }\zeta\in\mathcal{V}(\hat{\xi}).
\label{eq:thm5.1_b_1}%
\end{equation}

\end{enumerate}
\end{theorem}

\noindent{Proof.} \quad\text{ }\newline Statement b) is proved by reversing
the argument used to prove that \eqref{eq:thm5.1_b_1} $\Rightarrow$
\eqref{eq:thm5.1_a_1} - \eqref{eq:thm5.1_a_2}. We omit the details.
\hfill$\square$ \bigskip

\thanks{{\bf Acknowledgements}. N. Agram and B. \O ksendal are gratefully acknowledge the financial support provided by the Swedish Research Council grant (2020-04697) and the Norwegian Research Council grant (250768/F20), respectively.\\
The research of F. Proske was carried out with support of the Senter for internasjonalisering av utdanning (SIU), within the project Norway-Ukrainian  Cooperation in Mathematical Education, project number CPEA-LT-2016/10139.}

\end{document}